\documentclass[11pt,a4paper]{article}

\usepackage[utf8]{inputenc}
\usepackage[english]{babel}
\usepackage{amsmath}
\usepackage{amsthm}
\usepackage{amsfonts}
\usepackage{amssymb}
\usepackage{amscd}
\usepackage{tikz-cd}
\usepackage{bezier}
\usepackage{enumerate}
\usepackage{fancyhdr}
\usepackage{graphicx}
\usepackage{epstopdf}
\usepackage{booktabs}
\usepackage{rotating}
\usepackage{listings}
\usepackage[square, numbers, comma, sort&compress]{natbib} 
\usepackage{fancyhdr}
\usepackage{amsmath,amsfonts,amssymb,amscd,amsthm,xspace}
\usepackage[centerlast,small,sc]{caption}
\usepackage[pdfpagemode={UseOutlines},bookmarks=true,bookmarksopen=true,
   bookmarksopenlevel=0,bookmarksnumbered=true,hypertexnames=false,
   colorlinks,linkcolor={blue},citecolor={blue},urlcolor={red},
   pdfstartview={FitV},unicode,breaklinks=true]{hyperref}

\newtheorem{proposition}{Proposition}
\newtheorem*{proposition*}{Proposition}
\newtheorem*{theorem*}{Theorem}
\newtheorem{theorem}{Theorem}
\newtheorem{corollary}{Corollary}
\newtheorem{definition}{Definition}

\title{Positive topological entropy for Reeb flows on 3-dimensional Anosov contact manifolds}
\author{Marcelo R. R. Alves \\ \textsl{Institut Mittag-Leffler \& Universit\'{e} Paris-Sud}}

\begin{document}
\maketitle
\begin{center}
\abstract{\textbf{Abstract.} Let $(M, \xi)$ be a compact contact 3-manifold and assume that there exists a contact form $\alpha_0$ on $(M, \xi)$ whose Reeb flow is Anosov. We show this implies that every Reeb flow on $(M, \xi)$ has positive topological entropy, giving a positive answer to a question raised in \cite{A1}. Our
argument builds on previous work of the author \cite{A1} and recent work of Barthelm\'e
and Fenley \cite{BF}. \\ This result combined with the work of Foulon and Hasselblatt \cite{FH} is then used to obtain the first examples of hyperbolic contact 3-manifolds on which every Reeb flow has positive topological entropy.
}
\end{center}

\section{Introduction} \label{section1}

In this work we show that the existence of an Anosov Reeb flow on a compact contact 3-manifold $(M, \xi)$ has deep implications for the dynamics of all other Reeb flows on $(M, \xi)$. It implies that all Reeb flows on $(M, \xi)$ have positive topological entropy. It then follows from the works of Katok\cite{K,K2}, and Lima and Sarig \cite{LS} that every Reeb flow on $(M, \xi)$ has a compact invariant set on which the dynamics is conjugated to a subshift of finite type; i.e every Reeb flow on $(M, \xi)$ contains a ``horseshoe'' as a subsystem. Before giving the precise statements of our results we begin by recalling some necessary notions from contact geometry and dynamical systems.

The combination of the works of Klingenberg \cite{Klingenberg} and Manning \cite{Manning} implies that if a manifold $Q$ admits a Riemannian metric whose geodesic flow is an Anosov flow on the unit tangent bundle $T_1Q$, then for any Riemannian metric $g$ on $Q$ the associated geodesic flow has positive topological entropy. The results in the present paper can be seen as a generalization of this result to the world of contact 3-manifolds. It is not unreasonable to conjecture that an analogous result should hold for higher dimensional contact manifolds.

\subsection{Basic notions}

\subsubsection{Contact geometry and Reeb flows}
We first recall some basic definitions from contact geometry. A 1-form $\alpha$  on a $(2n+1)$-dimensional manifold $M$ called a \textit{contact form} if $ \alpha \wedge (d\alpha)^n $ is a volume form on $M$. The hyperplane $\xi= \ker \alpha$ is called the \textit{contact structure}. For us a \textit{contact manifold} will be a pair $(M,\xi)$ such that $\xi$ is the kernel of some contact form $\alpha$ on $M$ (these are usually called co-oriented contact manifolds in the literature). When $\alpha$ satisfies $\xi = \ker\alpha$, we will say that $\alpha$ is a contact form on $(M,\xi)$. On any contact manifold there always exist infinitely many different contact forms. Given a contact form $\alpha$, its
\textit{Reeb vector field} is the unique vector field $X_\alpha$ satisfying $\alpha(X_\alpha)=1$ and $i_{X_\alpha}d\alpha=0$. The Reeb flow $\phi_{X_{\alpha}}$ of $\alpha$ is the flow generated by the vector field $X_\alpha$. We will refer to the periodic orbits of $\phi_{X_{\alpha}}$ as \textit{Reeb orbits} of $\alpha$. The action $A(\gamma)$ of a Reeb orbit is defined by $A(\gamma):=\int_{\gamma} \alpha$. A contact form is called \textit{hypertight} if all its Reeb orbits are non-contractible.

To a non-degenerate Reeb orbit $\gamma$ of $\alpha$ and a trivialization $\Phi$ of the restriction $\gamma^*\xi$ of the contact distribution over $\gamma$, we associate an integer $\mu_{CZ}^{\Phi}(\gamma)$ called the Conley-Zehnder index of $\gamma$ with the trivialization $\Phi$; see for example \cite{B,Foliations} for the precise definition $\mu_{CZ}^{\Phi}(\gamma)$. Although $\mu_{CZ}^{\Phi}(\gamma)$ depends on the choice of the trivialization $\Phi$, its parity does not. We use this to define the parity of a Reeb orbit in the following way:
\begin{itemize}
  \item a non-degenerate Reeb orbit $\gamma$ is \textit{odd} if for all trivializations of $\gamma^*\xi$ its Conley-Zehnder index is odd,
  \item a non-degenerate Reeb orbit $\gamma$ is \textit{even} if for all trivializations of $\gamma^*\xi$ its Conley-Zehnder index is even.
\end{itemize}
It follows from our discussion that every non-degenerate Reeb orbit is either odd or even. 

\subsubsection{Dynamical systems} \label{dynamics}

Let $Y$ be a smooth 3-manifold and $X$ be a smooth vector field on $Y$ without singularities. We denote by $\phi_X$ the flow generated by $X$.
We endow $Y$ with an auxiliary Riemannian metric $g$ which generates a norm $|v|$ for every tangent vector $v$ in $M$. We will say that the flow $\phi_x$ is \textit{Anosov} when there exist line fields $E_S$ and $E_U$ in $Y$, and real numbers $\mu>0$ and $C>0$ such that:
\begin{itemize}
  \item at every point $y\in Y$ we have $T_y Y= \mathbb{R}X(y) \oplus E_S(y) \oplus E_U(y) $, where $\mathbb{R}X(y)$ is the 1-dimensional subspace of $T_yY$ generated by $X(y)$,
  \item the line fields $E_S$ and $E_U$ are left invariant by the flow $\phi_X$,
  \item $|D\phi_X^t(y)v| \geq Ae^{\mu t}|v|$ for every $y \in Y$, $v \in E_U(y)$ and $t\geq 0$,
  \item $|D\phi_X^t(y)v| \leq Ae^{-\mu t}|v|$ for every $y \in Y$, $v \in E_S(y)$ and $t\geq 0$.
\end{itemize}
Notice that it follows directly from these properties that all periodic orbits of an Anosov flow are hyperbolic.
An Anosov flow on a 3-manifold is called \textit{transversely orientable} when the line bundles $E_S$ and $E_U$ are trivial. 

Anosov flows occupy a very special place in the theory of dynamical systems. The first systematic study of this class of flows was done by 
Anosov \cite{Anosov} who showed that they are structurally stable. He also showed that if an Anosov flow preserves a smooth invariant measure that does not vanish on open sets then the flow is ergodic.

In this paper we prove a result about the positivity of the topological entropy for Reeb flows on certain contact manifolds. The topological entropy is a non-negative number associated to a dynamical system which measures the complexity of the orbit structure of the system. Positivity of the topological entropy means that the system possesses some type of exponential instability. We refer the reader to \cite{HK} for the definition of topological entropy. For flows on compact 3-manifolds that are generated by nowhere vanishing vector fields we have the following result of Katok \cite{K,K2}.
\begin{theorem*}[Katok\cite{K,K2}]
Let $X$ be a smooth vector field without singularities on a compact 3-manifold $M$, and let $\phi_X$ denote its flow. If $h_{top}(\phi_X) > 0$ then there exists a hyperbolic periodic orbit of $\phi_X$ whose unstable and stable manifolds intersect transversely. As a consequence there exists a compact set $K$ invariant by $\phi_X$, such that the dynamics of the restriction of $\phi_X$ to $K$ is conjugated to the suspension of a subshift of finite type.
\end{theorem*}
The recent work of Lima and Sarig \cite{LS} substantially improves this result. They showed that under the same hypothesis of Katok's theorem there exists a set $K$ invariant by the flow on which the dynamics can be described by a topological Markov chain and such that the topological entropy $h_{top}(\phi_X|_K)$ of the flow $\phi_X$ restricted to $K$ is equal to $h_{top}(\phi_X)$.

To estimate the topological entropy of a flow we use an idea present in \cite{A1}. Letting $M$ be a compact manifold and $X$ be a $C^2$ non-vanishing vector field on $M$, we define $\Lambda_X^T$ to be the set of free homotopy classes of loops in $M$ which contain a periodic orbit of the flow $\phi_X$ of $X$ with period $\leq T$. We denote by $N_X(T)$ the cardinality of $\Lambda_X^T$. We showed in \cite[Theorem 1]{A1} that if  $N_X(T)$ grows exponentially then the topological entropy of $\phi_X$ is positive.

\subsection{Main results}

Let $(M,\xi)$ be a compact contact 3-manifold. We will say that $(M,\xi)$ is an \textit{Anosov contact manifold} when there exists a contact form $\alpha_0$ on $(M,\xi)$ whose Reeb flow is Anosov. Anosov contact structures were studied in \cite{MP,Vaugon}. In \cite{MP} the authors establish the exponential growth of $S^1$-equivariant symplectic homology for symplectically fillable Anosov contact manifolds in all dimensions. In \cite{Vaugon} the author establish the exponential growth of cylindrical contact homology for Anosov contact 3-manifolds that admit a transversely orientable Anosov Reeb flow.

The classical examples of 3-dimensional Anosov contact 3-manifolds are the unit tangent bundles of higher genus surface endowed with the geodesic contact structure.\footnote{Positivity of topological entropy for Reeb flows on unit tangent bundles of higher genus surfaces endowed with the geodesic contact structure was obtained in \cite{MS} using a different approach based on volume growth of Lagrangian submanifolds; see also \cite{FS1,FS2}.} New interesting examples were recently constructed by Foulon and Hasselblatt in \cite{FH}.
Our main result is if $(M,\xi)$ is an Anosov contact 3-manifold then every Reeb flow on $(M,\xi)$ has positive topological entropy. We establish this result using a version of cylindrical contact homology. This is done via following theorem, which is a small modification of \cite[Theorem 2]{A1}.
\begin{theorem} \label{theorem1}
Let $\alpha_0$ be a hypertight contact form on a contact manifold $(M,\xi)$ and assume that the cylindrical contact homology of $\alpha_0$ has weak exponential homotopical growth with weight $a>0$. Then for every $C^k$ ($k\geq2$) contact form $\alpha$  on  $(M,\xi)$ the Reeb flow of $X_\alpha$ has positive topological entropy. More precisely, if $f_\alpha$ is the function such that $\alpha = f_\alpha \alpha_0$, then
\begin{equation}
h_{top}(\phi_{X_{\alpha}})\geq \frac{a}{\max f_\alpha}.
 \end{equation}
\end{theorem}

Our main theorem is:
\begin{theorem} \label{maintheorem}
Let $(M,\xi)$ be a compact 3-dimensional contact manifold and $\alpha_0$ be a contact form on $(M,\xi)$ such that its Reeb flow is
a transversely orientable Anosov flow. Then there exists a real number $a>0$ such that the cylindrical contact homology of $\alpha_0$
has weak exponential homotopical growth rate with weight $a>0$. It follows that if $\alpha$ is a contact form on $(M,\xi)$, and $f_\alpha$
is the function such that $f_\alpha \alpha_0= \alpha$, we have
\begin{equation}
h_{top}(X_\alpha)\geq \frac{a}{\max f_\alpha}.
\end{equation}
\end{theorem}
This theorem establishes the positive of topological entropy for all Reeb flows on an Anosov contact 3-manifold $(M,\xi)$ provided that $(M,\xi)$ admits a transversely orientable Anosov Reeb flow. The main tools in the proof of the theorem are the results of Fenley \cite{Fe} and, most importantly, the recent work of Barthelm\'e and Fenley \cite{BF}. The non-transversely orientable case is obtained as a corollary of the main theorem.
\begin{corollary} \label{corollary1}
Let $(M,\xi)$ be a compact contact 3-manifold and assume that there exists a contact form $\alpha_0$ on $(M,\xi)$ such that its Reeb flow is Anosov. Then every Reeb flow on $(M,\xi)$ has positive topological entropy.
\end{corollary}

Another corollary of Theorem \ref{maintheorem} is the first proof of existence of compact hyperbolic contact 3-manifolds\footnote{A compact hyperbolic contact 3-manifold is a compact contact 3-manifold $(M,\xi)$ such that $M$ is a compact hyperbolic 3-manifold.} on which every Reeb flow has positive topological entropy.
\begin{corollary} \label{corollary2}
There exists an infinity family $\{(M_n,\xi_n)\}_{n\in \mathbb{N}}$ of compact hyperbolic 3-manifolds such that:
\begin{itemize}
  \item $M_n$ and $M_m$ are not diffeomorphic if $n\neq m$,
  \item every Reeb flow on $(M_n,\xi_n)$ has positive topological entropy.
\end{itemize}
\end{corollary}
\textit{Proof of Corollary 2 assuming Theorem 1:}
In \cite{FH} Foulon and Hasselblatt showed that there exist infinitely many non-diffeomorphic hyperbolic 3-manifolds that admit 
contact forms whose Reeb flows are transversely orientable Anosov flows. This combined with Theorem \ref{maintheorem} implies the corollary.
\qed

\

A different approach to construct hyperbolic contact 3-manifolds on which every Reeb flow has positive topological entropy is the theme of a joint project of the author, Vincent Colin and Ko Honda, and is based on combining the methods of \cite{A2} and \cite{CH}.

\

\textbf{Structure of the paper: } In Section \ref{section2} we recall fundamental results about pseudoholomorphic curves in symplectizations and introduce the version of cylindrical contact homology that will be used in the paper. In Section \ref{section3} we prove Theorem \ref{theorem1}. Section \ref{section4} is devoted to the proofs of our main results, i.e Theorem \ref{maintheorem} and Corollary \ref{corollary1}.

\

\textbf{Acknowledgements:} I thank Fr\'ed\'eric Bourgeois and Anne Vaugon for many helpful conversations. This work was developed when the author was visiting the Mittag-Leffler Institute, which provided an extremely pleasant and rich environment for research.

\setcounter{theorem}{0}
\setcounter{corollary}{0}

\section{Pseudoholomorphic curves and cylindrical contact homology on free homotopy classes} \label{section2}

To define the cylindrical contact homology used in this paper we use pseudoholomorphic curves in symplectizations of contact manifolds and symplectic cobordisms. Pseudo-holomorphic curves were introduced in compact symplectic manifolds by Gromov in \cite{Gr} and adapted to symplectizations and symplectic cobordisms by Hofer \cite{H}. We refer the reader to \cite{CPT} as a general reference for pseudoholomorphic curves in symplectic cobordisms.

\subsection{Cylindrical almost complex structures}

Let $(M,\xi)$ be a contact manifold and $\alpha$ a contact form on $(M,\xi)$. The symplectization of $(M,\xi)$ is the product $\mathbb{R} \times M$ with the symplectic form $d(e^s \alpha)$ (where $s$ denotes the $\mathbb{R}$ coordinate in $\mathbb{R} \times M$). The 2-form $d\alpha$ restricts to a symplectic form on the vector bundle $\xi$ and it is well known that the set  $\mathfrak{j}(\alpha)$ of $d\alpha$-compatible almost complex structures on the symplectic vector bundle $\xi$ is non-empty and contractible. Notice that if $M$ is 3-dimensional, the set $\mathfrak{j}(\alpha)$ does not depend on the contact form $\alpha$ on $(M,\xi)$.

For $j \in \mathfrak{j}(\alpha)$ we can define an $\mathbb{R}$-invariant almost complex structure $J$ on $\mathbb{R} \times M$ by demanding that
\begin{equation}
J \partial_s = X_\alpha, \ \ J\mid_\xi = j.
\end{equation}
We will denote by $\mathcal{J}(\alpha)$ the set of almost complex structures in $\mathbb{R} \times M$ that are $\mathbb{R}$-invariant, $d(e^s\alpha)$-compatible and satisfy the equation (16) for some $j \in \mathfrak{j}(\alpha)$.

\subsection{Exact symplectic cobordisms with cylindrical ends}

An exact symplectic cobordism is, roughly, an exact symplectic manifold $(W,\varpi)$ that outside a compact subset is like the union of  cylindrical ends of symplectizations. We restrict our attention to exact symplectic cobordisms having only one positive end and one negative end.

Let $(W ,\varpi = d\kappa )$  be an exact symplectic manifold without boundary, and let $(M^+,\xi^+)$ and $(M^-,\xi^-)$ be contact manifolds with contact forms $\alpha^+$ and $\alpha^-$. We say that $(W ,\varpi = d\kappa )$ is an exact symplectic cobordism from  $\alpha^+$ to $\alpha^-$ if there exist subsets $W^-$, $W^+$ and $\widehat{W}$ of $W$ and diffeomorphisms $\Psi^+: W^+ \to [0,+\infty) \times M^+$ and $\Psi^-: W^- \to (-\infty,0] \times M^-$, such that:
\begin{eqnarray}
\widehat{W} \mbox{ is compact, } W= W^+ \cup \widehat{W} \cup W^- \mbox{ and } W^+ \cap W^- = \emptyset,
\end{eqnarray}
\begin{equation*}
(\Psi^+)^* (e^s\alpha^+) = \kappa \mbox{ and } (\Psi^-)^* (e^s\alpha^-) = \kappa.
\end{equation*}

In such a cobordism, we say that an almost complex structure $\overline{J}$ is cylindrical if
\begin{eqnarray}
\overline{J} \ \mbox{coincides with } \ J^+ \in \mathcal{J}(C^+ \alpha^+) \ \mbox{in the region} \ W^+, \\
\overline{J} \ \mbox{coincides with } \ J^- \in \mathcal{J}(C^- \alpha^-) \ \mbox{in the region} \ W^-, \\
\overline{J} \ \mbox{is compatible with} \ \varpi \ \mbox{in} \ \widehat{W},
\end{eqnarray}
where $C^+>0$ and $C^->0$ are constants.

For fixed $J^+ \in \mathcal{J}(C^+ \alpha^+)$ and $J^- \in \mathcal{J}(C^- \alpha^-)$, we denote by $\mathcal{J}(J^-,J^+)$ the
 set of cylindrical almost complex structures in $(\mathbb{R} \times M,\varpi)$ coinciding with $J^+$ on $W^+$ and $J^-$ on $W^-$.  It is well known that $\mathcal{J}(J^-,J^+)$ is non-empty and contractible. We will write $ \alpha^+ \succ_{ex} \alpha^-$ if there exists an exact symplectic cobordism from $ \alpha^+$ to $ \alpha^-$ as above. We remind the reader that $ \alpha^+ \succ_{ex} \alpha$ and $ \alpha \succ_{ex} \alpha^-$ implies $ \alpha^+ \succ_{ex} \alpha^-$; or in other words that the exact symplectic cobordism relation is transitive; see \cite{CPT} for a detailed discussion on symplectic cobordisms with cylindrical ends. Notice that a symplectization is a particular case of an exact symplectic cobordism.

\

\textbf{Remark:} We point out to the reader that in many references in the literature, a slightly different definition of cylindrical almost complex structures is used: instead of demanding that $\overline{J}$ satisfies equations (18) and (19), the stronger condition that $\overline{J}$ coincides with  $J^{\pm} \in \mathcal{J}( \alpha^{\pm}) \ \mbox{in the region} \ W^{\pm}$ is demanded. We need to consider this more relaxed definition of cylindrical almost complex structures when we study the cobordism maps of cylindrical contact homologies in subsection~\ref{sec3.3}.

\subsection{Splitting symplectic cobordisms}

Let $\alpha^+$, $\alpha$ and $\alpha^-$ be contact forms on $(M,\xi)$ such that $\alpha^+ \succ_{ex} \alpha$, $ \alpha \succ_{ex} \alpha^-$. For $\epsilon > 0$ sufficiently small, it is easy to see that one also has $\alpha^+ \succ_{ex} (1 + \epsilon)\alpha$ and $ (1-\epsilon)\alpha \succ_{ex} \alpha^-$. Then, for each $R>0$, we can to construct an exact symplectic form $\varpi_R = d\kappa_R$ on $W = \mathbb{R} \times M$ where:
\begin{eqnarray}
\kappa_R = e^{s-R-2}\alpha^+ \ \mbox{in} \ [R+2, + \infty) \times M, \\
\kappa_R = f(s)\alpha \ \mbox{in} \ [-R,R] \times M, \\
\kappa_R = e^{s+R+2}\alpha^- \ \mbox{in} \ (-\infty, -R-2] \times M,
\end{eqnarray}
and $f: [-R,R]  \to [1-\epsilon,1+\epsilon]$ satisfies $f(-R) = 1-\epsilon$, $f(R) = 1+\epsilon$ and $f'>0$.
In $ (\mathbb{R} \times M, \varpi_R)$ we consider a compatible cylindrical almost complex structure $\widetilde{J}_R$; but we demand an extra condition on $\widetilde{J}_R$:
\begin{eqnarray}
\widetilde{J}_R \ \mbox{coincides with} \ J \in \mathcal{J}(\alpha) \ \mbox{in} \ [-R,R] \times M.
\end{eqnarray}

Again we divide $W$ into regions:
$W^+ = [R+2, + \infty) \times M$, $W(\alpha^+,\alpha)= [ R, R+2] \times M$, $W(\alpha)=[-R, R] \times M$, $W(\alpha,\alpha^-)=[-R-2,-R] \times M$ and $W^-=(-\infty, -R-2] \times M$. The family of exact symplectic cobordisms with cylindrical almost complex structures $(\mathbb{R} \times M, \varpi_R, \widetilde{J}_R)$ is called a splitting family from $\alpha^+$ to $\alpha^-$ along $\alpha$.

\subsection{Pseudoholomorphic curves}

Let $(S,i)$ be a closed Riemann surface without boundary, $\Gamma \subset S$ be a finite set.
Let $\alpha$ be a contact form on $(M,\xi)$ and $J \in \mathcal{J}(\alpha)$. A finite energy pseudoholomorphic curve in the symplectization $(\mathbb{R} \times M,J)$ is a map $\widetilde{w}= (s,w):S \setminus \Gamma \to \mathbb{R} \times M$ that satisfies

\begin{eqnarray}
\overline{\partial}_J(\widetilde{w})= d\widetilde{w} \circ i - J \circ d\widetilde{w}=0
\end{eqnarray}
and
\begin{equation}
0<E(\widetilde{w})= \sup_{q \in \mathcal{E}} \int_{S \setminus \Gamma} \widetilde{w}^*d(q\alpha)
\end{equation}
where $\mathcal{E}= \{ q: \mathbb{R} \to [0,1]; q' \geq 0\}$. The quantity $E(\widetilde{w})$ is called the Hofer energy and was introduced in \cite{H}. The operator $\overline{\partial}_J$ above is called the Cauchy-Riemann operator for the almost complex structure $J$.

For an exact symplectic cobordism $(W,\varpi)$ from $\alpha^+$ to $\alpha^-$ as considered above, and $\overline{J} \in \mathcal{J}(J^-,J^+)$ a finite energy pseudoholomorphic curve is again a map $\widetilde{w}: S \setminus \Gamma \to W$ satisfying:
\begin{eqnarray}
d\widetilde{w} \circ i = \overline{J} \circ d\widetilde{w},
\end{eqnarray}
and
\begin{equation}
0<E_{\alpha^-}(\widetilde{w}) + E_c(\widetilde{w}) + E_{\alpha^+} (\widetilde{w}) < +\infty,
\end{equation}
where:
\\
$E_{\alpha^-}(\widetilde{w}) = \sup_{q \in \mathcal{E}} \int_{\widetilde{w}^{-1}(W^-))} \widetilde{w}^*d(q\alpha^-)$,
\\
$E_{\alpha^+}(\widetilde{w}) = \sup_{q \in \mathcal{E}} \int_{\widetilde{w}^{-1}(W^+)} \widetilde{w}^*d(q\alpha^+)$,
\\
$E_c(\widetilde{w}) = \int_{\widetilde{w}^{-1}(W(\alpha^+,\alpha^-))} \widetilde{w}^*\varpi$.
\\
These energies were also introduced in \cite{H}.

In splitting symplectic cobordisms we use a slightly modified version of energy. Instead of demanding $0<E_-(\widetilde{w}) + E_c(\widetilde{w}) + E_+ (\widetilde{w}) < +\infty$ we demand:
\begin{equation}
0<E_{\alpha^-}(\widetilde{w}) + E_{\alpha^-,\alpha}(\widetilde{w}) + E_{\alpha}(\widetilde{w})+ E_{\alpha,\alpha^+}(\widetilde{w})+ E_{\alpha^+} (\widetilde{w}) < +\infty
\end{equation}
where:
\\
$ E_{\alpha}(\widetilde{w})= \sup_{q \in \mathcal{E}} \int_{\widetilde{w}^{-1}W(\alpha)} \widetilde{w}^*d(q\alpha)$,
\\
$E_{\alpha^-,\alpha}(\widetilde{w})= \int_{\widetilde{w}^{-1}(W(\alpha,\alpha^-))} \widetilde{w}^*\varpi$, \
\\
$E_{\alpha,\alpha^+}(\widetilde{w}=\int_{\widetilde{w}^{-1}(W(\alpha^+,\alpha))} \widetilde{w}^*\varpi$,
\\
and $E_{\alpha^-}(\widetilde{w})$ and $E_{\alpha^+}(\widetilde{w})$ are as above.

\

The elements of the set $\Gamma \subset S$ are called punctures of the pseudoholomorphic curve $\widetilde{w}$. According to \cite{H, P1} punctures fall into two classes, positive and negative one according to the behaviour of $\widetilde{w}$ in the neighbourhood of the puncture.
Before presenting this classification we introduce some notation. Let $B_\delta (z)$ be the ball of radius $\delta$ centered at the puncture $z$, and denote by $\partial (B_\delta (z))$ its boundary. We can describe the types of punctures as follows:
\begin{itemize}
\item{$z \in \Gamma$ is called positive interior puncture if $z \in \Gamma$ and $\lim_{z' \to z} s(z') = +\infty$, and there exist a sequence $\delta_n \to 0$ and  Reeb orbit $\gamma^+$ of  $X_{\alpha^+}$, such that $w(\partial (B_{\delta_n} (z)))$ converges in $C^\infty $ to $\gamma^+$ as $n\to +\infty$}
\item{$z \in \Gamma$ is called negative interior puncture if $z \in \Gamma$ and $\lim_{z' \to z} s(z') = -\infty$, and there exist a sequence $\delta_n \to 0$ and  Reeb orbit $\gamma^-$ of  $X_{\alpha^-}$, such that $w(\partial (B_{\delta_n} (z)))$ converges in $C^\infty $ to $\gamma^-$ as $n\to +\infty$.}
\end{itemize}
The results in \cite{H} and \cite{P1}  imply that these are indeed the only possibilities we need to consider for the behaviour of $\widetilde{w}$ near punctures. Intuitively, we have that at the punctures, the pseudoholomorphic curve $\widetilde{w}$ detects Reeb orbits. When for a puncture $z$, there is a subsequence $\delta_n$ such that  $w(\partial(B_{\delta_n} (z)))$ converges to a Reeb orbit $\gamma$, we will say that $\widetilde{w}$ is asymptotic to this Reeb orbit $\gamma$ at the puncture $z$.

If a pseudoholomorphic curve  $\widetilde{w}$ is asymptotic to a non-degenerate Reeb orbit at a puncture $ z$, more can be said about its asymptotic behaviour in neighbourhoods of this puncture.
Take a neighbourhood $U \subset S$ of $z$ that admits a holomorphic chart $\psi_U : (U,z) \to (\mathbb{D},0)$. Using polar coordinates $(r,t) \in (0,+\infty) \times S^1$ we can write $x \in (\mathbb{D} \setminus 0)$ as $x = e^{-r}t $.
With this notation, it is shown in \cite{H,P1}, that if $z$ is a positive interior puncture at which $\widetilde{w}$ is asymptotic to a non-degenerate Reeb orbit $\gamma^+$ of $X_{\alpha^+}$, then $\widetilde{w} \circ \psi_U^{-1} (r,t) = (s(r,t),w(r,t))$ satisfies
\begin{itemize}
\item{$w^r(t)=w(r,t)$ converges in $C^{\infty}$ to a Reeb orbit $\gamma^+$ of $X_{\alpha^+}$, exponentially in $r$ and uniformly in $t$.}
\end{itemize}
Similarly, if $z$ is a negative interior puncture at which $\widetilde{w}$ is asymptotic to a non-degenerate Reeb orbit $\gamma^-$ of $X_{\alpha^-}$, then $\widetilde{w} \circ \psi_u^{-1} (r,t) = (s(r,t),w(r,t))$ satisfies
\begin{itemize}
\item{$w^r(t)=w(r,t)$ converges in $C^{\infty}$ to a Reeb orbit $\gamma^-$ of $-X_{\alpha^-}$ as  $r \to +\infty$, exponentially in $r$ and uniformly in $t$.}
\end{itemize}

\textit{Remark: The fact that the convergence of pseudoholomorphic curves near punctures to Reeb orbits is of exponential nature is a consequence of the asymptotic formula obtained in \cite{P1}. Such formulas are necessary for the Fredholm \cite{P2} theory that gives the dimension of the space of pseudoholomorphic curves with fixed asymptotic data.}

The discussion above can be summarised by saying that near punctures the finite energy pseudoholomorphic curves detect Reeb orbits. It is exactly this behavior that makes these objects useful for the study of dynamics of Reeb vector fields.

For us it will be important to consider the moduli spaces $\mathcal{M}(\gamma,\gamma'_1,...,\gamma'_m;J)$ of genus $0$ pseudoholomorphic curves, modulo biholomorphic reparametrisation,  with one positive puncture asymptotic to a non-degenerate Reeb orbit $\gamma$ and negative punctures asymptotic to non-degenerate Reeb orbits $\gamma'_1,...,\gamma'_m$. It is well known that the linearization $D\overline{\partial}_J$ of $\overline{\partial}_J$ at any element of $\mathcal{M}(\gamma,\gamma'_1,...,\gamma'_m;J)$ is a Fredholm map (we remark that this property is valid for more general moduli spaces of curves with prescribed asymptotic behaviour). One would like to conclude that the dimension of a connected component of $\mathcal{M}(\gamma,\gamma'_1,...,\gamma'_m;J)$ is given by the Fredholm index of an element of $\mathcal{M}(\gamma,\gamma'_1,...,\gamma'_m;J)$. However this is not always the case if the moduli space contain multiply covered pseudoholomorphic curves.

\textbf{Fact: As a consequence of the exactness of the symplectic cobordisms considered above we obtain that the energy $E(\widetilde{w})$ of $\widetilde{w}$ satisfies $E(\widetilde{w}) \leq 5A(\widetilde{w})$ where $ A(\widetilde{w})$ is the sum of the action of the Reeb orbits detected by the punctures of $\widetilde{w}$ counted with multiplicity; see \cite{CPT,Foliations}.}

\subsection{Cylindrical contact homology in special homotopy classes}\label{sec3.3}
We denote by $(M,\xi)$ a contact manifold endowed with a hypertight contact form $\alpha_0$.

Let $\Lambda_0$ denote the set of primitive free homotopy classes of $M$.  Let $\rho \in \Lambda$ be either an element of $\Lambda_0$, or a free homotopy class which contains only simple Reeb orbits of $\alpha_0$, and let $\mathcal{P}_\rho(\alpha_0)$ be the set of Reeb orbits of $\alpha_0$ that belong to $\rho$. Assume that all Reeb orbits in $\mathcal{P}_\rho(\alpha_0)$ are non-degenerate. By the work of Dragnev \cite{Dr}, we know that there exists a generic subset $\mathcal{J}^{\rho}_{reg}(\alpha_0)$ of $\mathcal{J}(\alpha_0)$ such that for all $J\in \mathcal{J}^{\rho}_{reg}(\alpha_0)$ we have:
\begin{itemize}
\item{for all Reeb orbits $\gamma_1,\gamma_2 \in \rho$, the moduli space of pseudoholomorphic cylinders $\mathcal{M}(\gamma_1,\gamma_2;J)$ is transverse, i.e. the linearized Cauchy-Riemann operator $D\overline{\partial}_J(\widetilde{w})$ is surjective for all $\widetilde{w} \in \mathcal{M}(\gamma_1,\gamma_2;J);$}
\item{for all Reeb orbits $\gamma_1,\gamma_2 \in \rho$, each connected component $\mathcal{L}$ of the moduli space $\mathcal{M}(\gamma_1,\gamma_2;J)$ is a manifold whose dimension is given by the Fredholm index of any element $\widetilde{w} \in \mathcal{L}$.}
\end{itemize}

For $\rho$ as above, we denote by $CH_{cyl}^{\rho}(\alpha_0)$ the $\mathbb{Q}$ vector-space generated by the set $\mathcal{P}_{\rho}(\alpha_0)$ of good Reeb orbits of $\alpha_0$ in $\rho$.
We then define, for $J \in \mathcal{J}^{\rho}_{reg}(\alpha_0)$:
\begin{equation}
d^{\rho}_J(\gamma)= cov(\gamma)\sum_{\gamma' \in \mathcal{P}_{\rho}(\alpha_0)} C^{\rho}(\gamma,\gamma';J)\gamma' = \sum_{\gamma' \in \mathcal{P}_{\rho}(\alpha_0)} C^{\rho}(\gamma,\gamma';J)\gamma'
\end{equation}
where $C^{\rho}(\gamma,\gamma';J)$ is the number of points of the moduli space $\widehat{\mathcal{M}}^1(\gamma,\gamma';J)$. The second equality follows from the fact that all Reeb orbits in $\rho$ are simple, which implies $cov(\gamma)=1$.

For $\alpha_0$ and $\rho$ as above and $J \in \mathcal{J}^{\rho}_{reg}(\alpha_0)$, the differential $d^{\rho}_{J}: CH_{cyl}^{\rho}(\alpha_0) \to CH_{cyl}^{\rho}(\alpha_0)$ is well-defined and satisfies $(d^{\rho}_{J})^2=0$. Thus, in this situation, we can define the cylindrical contact homology $C\mathbb{H}_{cyl}^{\rho,J}(\alpha_0)$. Once the transversality for $J$ has been achieved, and using coherent orientations constructed in \cite{BM}, the proof that $d^{\rho}_{J}$ is well-defined and that $(d^{\rho}_{J})^2=0$ is a combination of compactness and gluing, similar to the proof of the analogous result for Floer homology. The proof of these results is explained in \cite{A1,HMS,SFT}.
We will now state the results which we will use.
\begin{proposition} \label{proposition1} \cite{A1,HMS,SFT}
Let $(M,\xi)$ be a contact manifold with a hypertight contact form $\alpha_0$. Let $\rho \in \Lambda$ be either an element of $\Lambda_0$, or a free homotopy class which contains only simple Reeb orbits of $\alpha_0$. Assume that all Reeb orbits in $\mathcal{P}_\rho(\alpha_0)$ are non-degenerate and pick $J \in \mathcal{J}^{\rho}_{reg}(\alpha_0)$. Then, $d^{\rho}_{J}$ is well defined and $(d^{\rho}_{J})^2=0$. Under these conditions we define $C\mathbb{H}_{cyl}^{\rho}(\alpha_0)$ as the homology of the pair $(CH_{cyl}^{\rho}(\alpha_0), d^{\rho}_{J})$.
\end{proposition}
\textbf{Remark:} If we denote by $CH_{cyl|odd}^{\rho}(\alpha_0)$ and $CH_{cyl|even}^{\rho}(\alpha_0)$ the subspaces of $CH_{cyl}^{\rho}(\alpha_0)$ generated only by, respectively, odd and even Reeb orbits then it is clear that $CH_{cyl}^{\rho}(\alpha_0)= CH_{cyl|odd}^{\rho}(\alpha_0) \oplus CH_{cyl|even}^{\rho}(\alpha_0)$, and it is well known (see \cite{B,SFT}) that 
\begin{eqnarray}
  d^\rho_J(CH_{cyl|odd}^{\rho}(\alpha_0) ) \subset CH_{cyl|even}^{\rho}(\alpha_0)  \\
  d^\rho_J(CH_{cyl|even}^{\rho}(\alpha_0) ) \subset CH_{cyl|odd}^{\rho}(\alpha_0).
\end{eqnarray}

\begin{proposition} \label{proposition2}
Let $(M,\xi)$ be a contact manifold with a hypertight contact form $\alpha_0$. Let $\alpha^+ = C\alpha_0$ and $\alpha^- = c\alpha_0$ where $C>c>0$ are constants, and $\rho$ be either a primitive free homotopy class or a free homotopy class in which all Reeb orbits of $\alpha_0$ are simple. Assume that all Reeb orbits in $\mathcal{P}_\rho(\alpha_0)$ are non-degenerate.  Choose an almost complex structure $J\in \mathcal{J}^{\rho}_{reg}(\alpha_0)$, and set $J^+= J^- =J$. Let $(W=\mathbb{R} \times M,\varpi)$ be an exact symplectic cobordism from $C\alpha_0$ to $c\alpha_0$, and choose a regular almost complex structure $\widehat{J} \in \mathcal{J}^{\rho}_{reg}(J^-,J^+)$. Then, if there is a homotopy $(\mathbb{R}\times M,\varpi_t)$ of exact symplectic cobordisms from $C\alpha_0$ to $c\alpha_0$, with $\varpi_0 = \varpi$ and $\varpi_1= d(e^s\alpha_0)$, it follows that the map $\Phi^{\widehat{J}}:CH_{cyl}^{\rho,J}(\alpha_0) \to CH_{cyl}^{\rho,J}(\alpha_0) $ is chain homotopic to the identity.
\end{proposition}

\section{Growth of cylindrical contact homology and lower bounds for $h_{top}$} \label{section3}

In this section we use a modification of the notion of exponential homotopical growth of cylindrical contact homology which was introduced in \cite{A1}. This modification is needed in order to deal with the Anosov contact structures treated in the paper.

\begin{definition} \label{definitionhomotopy}
Let $(M,\xi)$ be a contact manifold and $\alpha_0$ be a hypertight contact form on  $(M,\xi)$. We denote by $\Lambda(M)$ the set of free homotopy classes of loops in $M$.
For $T>0$ let $\widetilde{\Lambda}_T(\alpha_0) \subset \Lambda(M)$ be the set which contains free homotopy classes $\rho$ that satisfy:
\begin{itemize}
\item all Reeb orbits of $X_{\alpha_0}$ in $\rho$ are simply covered, non-degenerate, and $C\mathbb{H}_{cyl}^{\rho}(\alpha_0) \neq 0$,
\item  there exists a finite subset $\{\gamma^{\rho}_1,...,\gamma^{\rho}_{k_\rho}\}\subset \mathcal{P}_{\rho}(\alpha_0)$ all of them with action $\leq T$ and non-zero rational numbers $q^{\rho}_1,...,q^{\rho}_{k_\rho}$ such that $d^{\rho}_{J}(\sum_{i=1}^{k_{\rho}} q^{\rho}_i \gamma^{\rho}_i)=0$ and $[\sum_{i=1}^{k_{\rho}} q^{\rho}_i \gamma^{\rho}_i] \in C\mathbb{H}_{cyl}^{\rho}(\alpha_0)$ is a non-zero homology class.
\end{itemize}
We define $N^{cyl}_T(\alpha_0) := \# \widetilde{\Lambda}_T(\alpha_0)$.
\end{definition}
The second condition, which is the crucial one, means that the there is a non-vanishing homology class in $C\mathbb{H}_{cyl}^{\rho}(\alpha_0)$ which can be represented as a sum of Reeb orbits of $\alpha_0$ with the action $\leq T$.
\begin{definition} \label{definition2}
We say that the cylindrical contact homology of $\alpha_0$ of $(M,\alpha_0)$ has weak exponential homotopical growth with weight $a>0$ if there exist a number $b$ and a sequence $T_n \to +\infty$, such that $N^{cyl}_{T_n}(\alpha_0) \geq e^{aT_n + b}$ for all $T_n$.
\end{definition}

The following theorem is a strengthening of \cite[Theorem 2]{A1}:
\begin{theorem} \label{theorem1'}
Let $\alpha_0$ be a hypertight contact form on a contact manifold $(M,\xi)$ and assume that the cylindrical contact homology of $\alpha_0$ has weak exponential homotopical growth with weight $a>0$. Then for every $C^k$ ($k\geq2$) contact form $\alpha$  on  $(M,\xi)$ the Reeb flow of $X_\alpha$ has positive topological entropy. More precisely, if $f_\alpha$ is the function such that $\alpha = f_\alpha \alpha_0$, then
\begin{equation}
h_{top}(\phi_{X_{\alpha}})\geq \frac{a}{\max f_\alpha}.
 \end{equation}
\end{theorem}
\textit{Proof:} The proof is almost identical to the one of \cite[Theorem 2]{A1}.
We write $E=\max f_{\alpha}$.

\textbf{Step 1:} \\
 We assume first that $\alpha$ is non-degenerate and $C^{\infty}$.
For every $\epsilon>0$ we can construct an exact symplectic cobordism from $(E+\epsilon) \alpha_0$ to $\alpha$. Analogously, for $e>0$ small enough, it is possible to construct an exact symplectic cobordism from $\alpha$ to $e\alpha_0$.

Let $\rho \in \widetilde{\Lambda}_T(\alpha_0)$.
Using these cobordisms, we can construct a splitting family $(\mathbb{R} \times M, \varpi_R, J_R)$ from $(E+\epsilon) \alpha_0$ to $e\alpha_0$, along $\alpha$, such that for every $R>0$, $(\mathbb{R} \times M, \varpi_R, J_R)$ is homotopic to the symplectization of $\alpha_0$ and $J_R$ is regular. To do this, we first fix almost complex structures $J_0 \in \mathcal{J}^{\rho}_{reg}(\alpha_0)$ and $J \in \mathcal{J}(\alpha)$, and pick an almost complex structure $J'_R$ compatible with $\varpi_R$ and that coincides with $J_0$ on the positive and negative ends of the cobordism, and with $J$ on $[-R,R] \times M$. It follows from the methods of \cite{Bourgeoishomotopygroup} that we can make a perturbation of $J'_R$ supported in $(R,R+1)\times M$ to obtain $J_R$ compatible with $\varpi_R$ that belongs to $\mathcal{J}^{\rho}_{reg}(J_0,J_0)$ and coincides with $J$ on $[-R,R] \times M$.

Since $J_R \in \mathcal{J}^{\rho}_{reg}(J_0,J_0)$, the exact symplectic cobordism  $(\mathbb{R} \times M, \varpi_R, J_R)$ induces a map $\Psi_{J_R}: C\mathbb{H}_{cyl}^{\rho}(\alpha_0) \to C\mathbb{H}_{cyl}^{\rho}(\alpha_0)$. Because $(\mathbb{R} \times M, \varpi_R, J_R)$ is homotopic to the symplectization of $\alpha_0$, we know from Proposition \ref{proposition2} that $\Psi_{J_R}$ is the identity.

For our $\rho \in \widetilde{\Lambda}_T(\alpha_0)$, let $\{\gamma^{\rho}_1,...,\gamma^{\rho}_{k_\rho}\}\subset \mathcal{P}_{\rho}(\alpha_0)$ and $q^{\rho}_1,...,q^{\rho}_{k_\rho}$ be, respectively, the Reeb orbits with action $\leq T$ and rational numbers given in Definition \ref{definition2}, such that $d^{\rho}_{J}(\sum_{i=1}^{k_{\rho}} q^{\rho}_i \gamma^{\rho}_i)=0$ and $[\sum_{i=1}^{k_{\rho}} q^{\rho}_i \gamma^{\rho}_i] \in C\mathbb{H}_{cyl}^{\rho}(\alpha_0)$ is non-zero. Since $\Psi_{J_R}: C\mathbb{H}_{cyl}^{\rho}(\alpha_0) \to C\mathbb{H}_{cyl}^{\rho}(\alpha_0)$ is the identity, it follows that $\Psi_{J_R}([\sum_{i=1}^{k_{\rho}} q^{\rho}_i \gamma^{\rho}_i]) \neq 0$.
We then conclude that for every $R$ there exists some $i_R \in \{1,...,k_{\rho}\}$ and a finite energy pseudoholomorphic cylinder $\widetilde{w}$ in $(\mathbb{R} \times M,J_R)$ positively asymptotic to $\gamma^\rho_{i_R}$ and negatively asymptotic to an orbit in $\mathcal{P}_\rho(\alpha_0)$.

\

\textbf{Step 2:} \\
Let $\rho \in \widetilde{\Lambda}_T(\alpha_0) $, $R_n \to +\infty$ be a strictly increasing sequence. For each $R_n$, let $i_{R_n} \in \{1,...,k_{\rho} $ be the number and $\widetilde{w}_n : (S^1 \times \mathbb{R},i) \to (\mathbb{R} \times M, J_{R_n}) $ be the pseudoholomorphic cylinders with one positive puncture asymptotic to an orbit in $\gamma^{\rho}_{i_{R_n}}$ and one negative puncture asymptotic to an orbit in $\mathcal{P}_\rho(\alpha_0)$, which were obtained in Step 1. Notice that, because the action of $\gamma^{\rho}_{i_{R_n}}$ as Reeb orbit of $(E+\epsilon)\alpha_0$ is $\leq (E+\epsilon)T$, we know that the energy of $\widetilde{w}_n $ is uniformly bounded.

Therefore we can apply the SFT compactness theorem of \cite{CPT} to obtain a subsequence of $\widetilde{w}_n $ which converges to a pseudoholomorphic buiding $\widetilde{w}$. Notice that in order to apply the SFT compactness theorem we need to use the non-degeneracy of $\alpha$. Moreover we can give a very precise description of the building.

Let $\widetilde{w}^k$ for $k\in \{1,...,m\}$ be the levels of the pseudoholomorphic building $\widetilde{w}$. Because the topology of our curve $\widetilde{w}$ does not change after breaking we have the following picture:
\begin{itemize}
\item{the upper level $\widetilde{w}^1$ is composed of one connected pseudoholomorphic curve, which has one positive puncture asymptotic to $\gamma^{\rho}_{i_{R_n}}$, and several negative punctures. All of the negative punctures detect contractible orbits, except one that detects a Reeb orbit $\gamma_1$ which is also in $\rho$.}
\item{on every other level $\widetilde{w}^k$ there is a special pseudoholomorphic curve which has one positive puncture asymptotic to a Reeb orbit $\gamma_{k-1}$ in $\rho$, and at least one but possibly several negative punctures. Of the negative punctures there is one that is asymptotic to an orbit $\gamma_k$ in $\rho$, while all the others detect contractible Reeb orbits.}
\end{itemize}

Because of the splitting behavior of the cobordisms $(\mathbb{R} \times M, J_{R_n})$ it is clear that there exists a $k_0$, such that the level $\widetilde{w}^{k_0}$ is in an exact symplectic cobordism from $(E+\epsilon)\alpha_0$ to $\alpha$. This implies that the special orbit $\gamma_{k_0}$ is a Reeb orbit of $X_\alpha$ in the homotopy class $\rho$.

As the action of $\gamma^{\rho}_{i_{R_n}}$ as Reeb orbit of $(E+\epsilon)\alpha_0$ is $\leq (E+\epsilon)T$, we conclude that $A(\gamma_0) \leq (E+\epsilon)T$. This implies that all the other orbits appearing as punctures of the building $\widetilde{w}$ have action smaller than $(E+\epsilon)T$, and in particular that $\gamma_{k_0}$ has action smaller than $(E+\epsilon)T$.

As we can do the construction above for any $\epsilon>0$ we can obtain a sequence of Reeb orbits $\gamma^{\rho}_{j}$ which are all in $\rho$ and such that $A(\gamma^{\rho}_{j}) \leq (E+\frac{1}{j})T$. Using the Arzela-Ascoli Theorem one can extract a convergent subsequence of $\gamma_j^\rho$. Its limit $\gamma_\rho$ is clearly a Reeb orbit of $\alpha$ in the free homotopy class $\rho$ and with action $\leq ET$.

\

\textbf{Step 3:} Estimating the number $N_{X_{\alpha}}(T)$ of free homotopy classes of $M$ that contain periodic orbits of $X_\alpha$ with period $\leq T$; see the Subsection \ref{dynamics} for the definition $N_{X_{\alpha}}(T)$.

From Step 2, we know that if $\rho \in \widetilde{\Lambda}_T(\alpha_0)$ then there is a Reeb orbit $\gamma_\rho$ of the Reeb flow of $X_\alpha$ with $A(\gamma_\rho)\leq ET$. Recalling that the period and the action of a Reeb orbit coincide we obtain that $N_{X_{\alpha}}(T) \geq \#\widetilde{\Lambda}_{\frac{T}{E}}(\alpha_0)$. Under the hypothesis of the theorem there exists a sequence $T_n \to +\infty$ such that $\#\widetilde{\Lambda}_{\frac{T}{E}}(\alpha_0) \geq e^{\frac{aT_n}{E} + b}$ for all $T_n$. We then conclude that
\begin{equation}
N_{X_{\alpha}}(T_n) \geq e^{\frac{aT_n}{E} + b}
\end{equation}
for all elements of the sequence $T_n$.
Applying \cite[Theorem 1]{A1} we obtain $h_{top}(\phi_{X_{\alpha}}) \geq \frac{a}{E}$. This proves the theorem in the case that $\alpha$ is $C^{\infty}$ and non-degenerate.

\

\textbf{Step 4:} Passing to the case of a general $C^{k\geq2}$ contact form $\alpha$ (the case where $\alpha$ is degenerate is included here).

Let $\alpha_i$ be a sequence of non-degenerate smooth contact forms converging in the $C^k$-topology to a contact form $\alpha$ which is $C^k $ ($k\geq 2$) and possibly degenerate. For every $\epsilon >0$ there is $i_0$ such that for $i>i_0$ there exists an exact symplectic cobordism from $(E+\epsilon)\alpha_0 $ to $\alpha_i$.

Fixing then a homotopy class $ \rho \in \widetilde{\Lambda}_T(\alpha_0)$ we know, by the previous steps, that there exists a Reeb orbit $\gamma_\rho(i)$ of $\alpha_i$ in the homotopy class $\rho$ with action smaller than $(E+\epsilon)T$. By applying the Arzela-Ascoli theorem to $\gamma_\rho(i)$, we obtain a subsequence which converges to a Reeb orbit $\gamma_{\epsilon,\rho}$ of $X_\alpha$ with $A(\gamma_{\epsilon,\rho})\leq (E+\epsilon)T$.

Because $\epsilon>0$ above can be taken arbitrarily close to $0$ we can actually obtain a sequence $\gamma_{j,\rho}$ of Reeb orbits of $X_\alpha$ whose homotopy class is $\rho$ such that the actions $A(\gamma_{j,\rho})$ converges to a number $\leq ET$. Again applying Arzela-Ascoli theorem, we obtain that the sequence $\gamma_{j,\rho}$ has a convergent subsequence, which converges to an orbit $\gamma_{\rho}$ satisfying $A(\gamma_{\rho}) \leq ET$.

Reasoning as in Step 3 above, we conclude that $N_{X_{\alpha}}(T_n) \geq e^{\frac{aT_n}{E} + b}$ for all elements of the sequence $T_n \to +\infty$. Applying Theorem \ref{theorem1'} we obtain the desired estimate for the topological entropy. This finishes the proof of the theorem.
\qed

\section{Weak exponential homotopical growth of cylindrical contact homology for Anosov Reeb flows} \label{section4}

In this section we prove Theorem \ref{maintheorem}.
We start by recalling some facts about Anosov flows on compact 3-dimensional manifolds. It follows from Novikov's theorem that all the
periodic orbits of a 3-dimensional Reeb flows are non contractible (see \cite{Fe}). In \cite{Fe} Fenley showed that if a free homotopy class of a 3-manifold contains a simple periodic orbit of a transversely orientable\footnote{Recall that an Anosov flow in a 3-manifold
is transversely orientable if the strong stable and unstable bundles are trivial subbundles of the tangent bundle of the manifold.}
Anosov flow, then the free homotopy class is primitive. Combining these two results we obtain the following
\begin{proposition}[Anosov\cite{Anosov},Fenley\cite{Fe}] \label{propositionAF}
Let $(M,\xi)$ be a compact 3-dimensional contact manifold and $\alpha_0$ be a contact form on $(M,\xi)$ such that its Reeb flow is a transversely orientable Anosov flow.
Then $\alpha_0$ is hypertight and all simple Reeb orbits of $\alpha_0$ are contained in primitive free homotopy classes.
\end{proposition}

Crucial for the proof of Theorem \ref{maintheorem} is the recent work of Barthelm\'e and Fenley \cite{BF}. They proved that the number of free
homotopy classes in a compact 3-manifold containing periodic orbits with action $\leq T$ of an Anosov flow grow exponentially with $T$.
The application of their result to transversely orientable Anosov Reeb flows on compact contact 3-manifolds gives us
\begin{proposition} \label{mainproposition}
Let $(M,\xi)$ be a compact 3-dimensional contact manifold and $\alpha_0$ be a contact form on $(M,\xi)$ such that its Reeb flow is
a transversely orientable Anosov flow. Let $\widetilde{\Lambda}_0^T(\alpha_0)$ be the set of primitive free homotopy classes of $M$
that contain a Reeb orbit of $\alpha_0$ with action $\leq T$. Then there exists a monotone sequence $T_n \to +\infty$ and real numbers
$a>0$ and $b$ such that
\begin{equation}
\#\Lambda_0^{T_n}(\alpha_0) > e^{aT_n + b}. \end{equation}
\end{proposition}
\textit{Proof:}
We first define the set $\Lambda^{T}(\alpha_0)$ as the set of free homotopy classes in $M$ that contain some Reeb orbit of $\alpha_0$ with action $\leq T$.

As the Reeb vector field $X_{\alpha_0}$ does not have singularities and $M$ is compact, we know that there exists a constant $c>0$ such that all Reeb orbits of $\alpha_0$ have action $\geq c$.

Let $\rho \in \Lambda^{T}(\alpha_0)$. Then it follows from Proposition \ref{propositionAF} that we have the following dichotomy:
\begin{itemize}
\item $\rho$ is primitive and therefore belongs to $\Lambda_0^T(\alpha_0)$,
\item or $\rho$ contains only multiply covered Reeb orbits, which are all multiple covers of a simple Reeb orbits that belong to
      a primitive homotopy class $\rho' \in \Lambda_0^T(\alpha_0)$.                                                               \end{itemize}

Let $\rho$ be an element of $\Lambda_0^T(\alpha_0)$. Given a positive integer $n$ we will denote by $\rho^n$ the free homotopy class
that contains $k$-covers of closed curves in $\rho$. With this notation we claim that if $\rho^k \in \Lambda^{T}(\alpha_0)$ then we must have
$k\leq \lfloor\frac{T}{c}\rfloor$. To see that this is indeed the case, note that if $\rho^k \in \Lambda^{T}(\alpha_0)$ then it follows from the above mentioned dichotomy that every Reeb orbit $\gamma$ in $\rho^k$ is the $k$-cover of some Reeb orbit $\widehat{\gamma}$ in $\rho$. As   $\gamma=\widehat{\gamma}^k$ we obtain that $A(\gamma)= A(\widehat{\gamma}^k)\geq ck$. We have showed that all Reeb orbits in $\rho^k$ have       action $\geq ck$, therefore as $\rho^k \in \Lambda^{T}(\alpha_0)$ we must have $ck \leq T$ which implies our claim that $k\leq \lfloor\frac{T}{c}\rfloor$.

Our discussion so far implies the following situation:
\begin{itemize}
  \item if $\rho \in \Lambda^{T}(\alpha_0)$ then there exist a free homotopy class $\varrho \in \Lambda_0^T(\alpha_0)$ and
  an integer $k \in \{1,...,\lfloor\frac{T}{c}\rfloor \}$ such that $\varrho^k = \rho$.
\end{itemize}
This implies that
\begin{equation} \label{estimate1}
\lfloor\frac{T}{c}\rfloor \#\Lambda_0^{T}(\alpha_0) \geq \# \Lambda^{T}(\alpha_0).
\end{equation}

We now invoke \cite[Theorem A]{BF}. It implies that there exist a monotone sequence $T'_n \to +\infty$ and real numbers $\widetilde{a}>0$ and $b$
such that
\begin{equation} \label{estimate2}
\# \Lambda^{T'_n}(\alpha_0) \geq e^{\widetilde{a}T'_n + b}.
\end{equation}
Combining equations \ref{estimate1} and \ref{estimate2} we obtain that
\begin{equation} \label{estimate3}
\#\Lambda_0^{T'_n}(\alpha_0) > \frac{e^{\widetilde{a}T'_n + b}}{\lfloor\frac{T'_n}{c}\rfloor}.
\end{equation}

As $\lfloor\frac{T}{c}\rfloor$ grows linearly with $T$, it follows from equation \ref{estimate3} that for any if we pick a positive real number $a<\widetilde{a}$, then there exists a monotone sequence $T_n \to +\infty$ such that
\begin{equation}\label{estimate4}
\#\Lambda_0^{T_n}(\alpha_0) \geq e^{aT_n + b}.
\end{equation}
\qed

We are now ready to prove Theorem \ref{maintheorem}.
\begin{theorem} \label{maintheorem'}
Let $(M,\xi)$ be a compact 3-dimensional contact manifold and $\alpha_0$ be a contact form on $(M,\xi)$ such that its Reeb flow is
a transversely orientable Anosov flow. Then there exists a real number $a>0$ such that the cylindrical contact homology of $\alpha_0$
has weak exponential homotopical growth rate with weight $a>0$. It follows that if $\alpha$ is a contact form on $(M,\xi)$, and $f_\alpha$
is the function such that $f_\alpha \alpha_0= \alpha$, we have
\begin{equation}
h_{top}(X_\alpha)\geq \frac{a}{\max f_\alpha}.
\end{equation}
\end{theorem}
\textit{Proof:} Let $a>0$ be as in statement of Proposition \ref{mainproposition}.

\textbf{Step 1:} \\
We first show that for every $\rho \in \Lambda_0^{T}(\alpha_0)$ we have $C\mathbb{H}_{cyl}^{\rho}(\alpha_0) \neq 0$.
To show that we first notice that since $\rho$ is primitive and $\alpha_0$ are non-degenerate\footnote{As the Reeb flow of $\alpha_0$ is Anosov all its Reeb orbits are hyperbolic and therefore non-degenerate.},
we can pick an almost complex structure $J \in \mathcal{J}^\rho(\alpha_0)$ for which differential $d^\rho_J$ is well-defined.
This implies that $C\mathbb{H}_{cyl}^{\rho}(\alpha_0)$.

Since the Reeb flow of $\alpha_0$ is Anosov all its Reeb orbits are hyperbolic. Moreover, the fact that Reeb flow of $\alpha_0$ is a transversely orientable Anosov flow implies that for any Reeb orbit $\gamma$ of $\alpha_0$ its stable and unstable bundles are trivial line bundles over $\gamma$. This implies that all the Reeb orbits of $\alpha_0$ are hyperbolic and even; see for example \cite{Vaugon}.

It follows from the Remark following Proposition \ref{proposition1}, that for any even Reeb orbit $\gamma \in \rho$ the differential $d^\rho_J(\gamma)$ is a sum of odd Reeb orbits in $\rho$.
As there $\alpha_0$ has no odd Reeb orbits we conclude that the differential $d^\rho_J$ is the zero map. This implies that the chain complex $CH_{cyl}^{\rho}(\alpha_0)$ is isomorphic  to the homology $C\mathbb{H}_{cyl}^{\rho}(\alpha_0)$ for all $\rho \in \Lambda_0^{T}(\alpha_0)$. Since for all such $\rho$ we have $CH_{cyl}^{\rho}(\alpha_0)\neq 0$, we conclude that for all $\rho \in \Lambda_0^{T}(\alpha_0)$ the homology
$C\mathbb{H}_{cyl}^{\rho}(\alpha_0) \neq 0$.

\

\textbf{Step 2:} \\
In Step 1 we showed that for all $\rho \in \Lambda_0^{T}(\alpha_0)$, the differential $d^\rho_J:CH_{cyl}^{\rho}(\alpha_0) \to CH_{cyl}^{\rho}(\alpha_0)$ is the zero map. This implies that for each Reeb orbit $\gamma \in \rho$ is a closed and non-exact element
of the chain complex $(CH_{cyl}^{\rho}(\alpha_0),d^\rho_J)$. As a conclusion we have that each Reeb orbit $\gamma \in \rho$ is the representative of a non-zero homology class $[\gamma] \in CH_{cyl}^{\rho}(\alpha_0)$.

\

\textbf{Step 3:} We will show that $\Lambda_0^{T}(\alpha_0) \subset \widetilde{\Lambda}_T(\alpha_0)$, with $\widetilde{\Lambda}_T(\alpha_0)$
as defined in Definition \ref{definitionhomotopy}. \\
Let $\rho \in \Lambda_0^{T}(\alpha_0)$, and let $\gamma_\rho$ be a Reeb orbit of $\alpha_0$ in $\rho$ with action $\leq T$.
It follows from the previous steps that $C\mathbb{H}_{cyl}^{\rho}(\alpha_0) \neq 0$, $d^\rho_J(\gamma_\rho)=0$ and the homology class
$[\gamma_\rho]$ is non-zero in $C\mathbb{H}_{cyl}^{\rho}(\alpha_0)$. It follows then that $\rho \in \widetilde{\Lambda}_T(\alpha_0)$ as claimed.

\textbf{Step 4:} \\
As $\Lambda_0^{T}(\alpha_0) \subset \widetilde{\Lambda}_T(\alpha_0)$, it follows that for the monotone sequence $T_n \to +\infty$ and the numbers
$a>0$ and $b$ given in the statement of Proposition \ref{mainproposition} we have:
\begin{equation}\label{estimatefinal}
N^{cyl}_{T_n}(\alpha_0) = \#\widetilde{\Lambda}_{T_n}(\alpha_0) \geq \#\Lambda_0^{T_n}(\alpha_0) \geq e^{aT_n + b}, \end{equation}
which finishes the proof of the first statement of theorem. The second statement of the theorem follows from combining the first with Theorem \ref{maintheorem'}.
\qed

 Theorem \ref{maintheorem'} we obtain the following:
\begin{corollary} \label{corollary1'}
Let $(M,\xi)$ be a compact contact 3-manifold and assume that there exists a contact form $\alpha_0$ on $(M,\xi)$ such that its Reeb flow is Anosov. Then every Reeb flow on $(M,\xi)$ has positive topological entropy.
\end{corollary}
\textit{Proof:} \\
If the Reeb flow of  $\alpha_0$  is a transversely orientable Anosov flow then we showed in Theorem \ref{maintheorem'} that every
Reeb flow on $(M,\xi)$ has positive topological entropy.

We will thus assume from now on that this is not the case. Then we know \cite{Fenley1} that there exists a double covering $\pi_2 :\widetilde{M} \to M$ such that the flow of the pullback vector field $\pi_2^* X_{\alpha_0}$ is a transversely orientable Anosov flow.
It is clear that $\widetilde{\xi}= \pi_2^* \xi$ is a
contact structure on $\widetilde{M}$ and that $\widetilde{\alpha}_0 := \pi_2^* \alpha_0 $ is a contact form on $(\widetilde{M},\widetilde{\xi})$. The Reeb vector field $X_{\widetilde{\alpha}_0}$ of $\widetilde{\alpha}_0$ coincides with $\pi_2^* X_{\alpha_0}$. We thus conclude that the Reeb flow of $\widetilde{\alpha}_0$ is a transversely orientable Anosov flow, and it follows from
\ref{maintheorem'} that every Reeb flow on $(\widetilde{M},\widetilde{\xi})$ has positive topological entropy.

Taking $\alpha$ to be a contact form on $(M,\xi)$, we define $\widetilde{\alpha}:= \pi_2^* \alpha$. It is clear that $\widetilde{\alpha}$ is a
contact form on $(\widetilde{M},\widetilde{\xi})$, which implies that its Reeb flow has positive topological entropy. It follows directly that
the Reeb flow of $\alpha$ must also have positive topological entropy.
\qed

\end{document}